\documentclass[11pt,a4paper]{article}
\usepackage{amsfonts,amsmath,amstext,amssymb}
\usepackage{dsfont}
\usepackage{theorem}    %para cambiar el estilo de los teoremas
\usepackage{a4wide}
\usepackage{verbatim}

\newtheorem{lema}{Lemma}
\newtheorem{prop}[lema]{\bf Proposition}
\newtheorem{teor}[lema]{\bf Theorem}
\newtheorem{coro}[lema]{\bf Corollary}

\newtheorem{rema}[lema]{\bf Remark}

\theorembodyfont{\rmfamily}

\title{Compact maximal hypersurfaces in stably causal spacetimes}

\author{\footnote{The authors are partially supported by the
Spanish MEC-FEDER Grant MTM2010-18099. }
Rafael M. Rubio and  Juan J. Salamanca\\[6mm]
 Departamento de Matem\'aticas, Campus de Rabanales, \\[0.5mm]
Universidad de C\'ordoba, 14071 C\'ordoba, Spain,\\[0.5mm]
E-mails\textup{:\texttt{\;rmrubio@uco.es},\,\,\texttt{jjsalamanca@uco.es}}}

\date{}

\begin{document}

\maketitle

\thispagestyle{empty}

\begin{abstract}
Several uniqueness results on compact maximal hypersurfaces in a wide class of stably causal spacetimes are given. They are obtained from the study of a distinguished function on the maximal hypersurface, under suitable natural first order conditions of the spacetime. As a consequence several applications to Geometric Analysis are given.
\vspace*{2mm}

\end{abstract}


\noindent {{\bf{PACS 2010:}} 02.40.Ky, 04.20.Cv.} 

\vspace*{-2mm}

\section{Introduction}

A continous function on a spacetime, which is strictly increasing on any future directed causal curve is called a time function. In 1968,  Hawking proved that any stably causal spacetime admits a time function \cite{Haw}. Conversely, Hawking also proved that a spacetime is stably causal if it admits a smooth function with a everywhere timelike gradient. More recently, Bernal and Sanchez in \cite{BS} define a temporal function as   a smooth function on the spacetime with (past-directed) timelike gradient on all  the spacetime and prove that any spacetime which admits a time function also admits a temporal function. Thus, as the authors highlight, this result combined with Hawking's ones, ensures, on one hand, that any stable causal spacetime admits a temporal function and, on the other, that any spacetime which admits a time function is stably causal.

On the other hand, there exist a serie of folk questions refered to differentiability of time functions and
Caucly hypersurfaces and smooth orthogonal splitting of globally hyperbolic spacetimes. The principal Theorem in \cite{BS} assures that an $(n+1)$-dimensional globally hyperbolic 
spacetime $M$ is isometric to the smooth product manifold 
\begin{equation} \label{ghg}
\mathbb{R}\times \mathcal{S} \, , \quad \overline{g}= -\beta d\mathcal{T}^2+\hat{g} 
\end{equation} where $\mathcal{S}$ is a smooth spacelike Cauchy hypersurface, $\mathcal{T}: \mathbb{R}\times \mathcal{S}\rightarrow \mathbb{R}$
is the natural projection, $\beta: \mathbb{R}\times \mathcal{S}\rightarrow (0,\infty)$ a smooth function and $\hat{g}$ a $2$-covariant
symmetric tensor field on $\mathbb{R}\times \mathcal{S}$, satisfying:

i) $\nabla{\mathcal{T}}$ is timelike and past-pointing on all $M$.

ii) Each hypersurface $\mathcal{S}_{\mathcal{T}}$ at constant ${\mathcal{T}}$ is a Cauchy hypersurface, and the restriction $\hat{g}_{_{\mathcal{T}}}$ of $\hat{g}$ to such a $\mathcal{S}_{\mathcal{T}}$ is a Riemannian metric.

iii) The radical of $\hat{g}$ at each $w\in \mathbb{R}\times \mathcal{S}$ is $Span \nabla {\mathcal{T}}$ at $w$.

\vspace{3mm}

If we denote $\pi_{_\mathcal{S}}$ the canonical projection on the second factor, which is a diffeomorphism restricted to each 
level hypersurface $\mathcal{S}_{\mathcal{T}}$, via $\pi_{_\mathcal{S}}$, we can give an one-parameter family of Riemannian metric $g_{_\mathcal{T}}$ on the differentiable manifold $\mathcal{S}$.

Inspired in all these facts, we consider a natural family of stably causal spacetimes which allows certain accouracy to the study of maximal hypersurfaces.

\vspace{3mm}

Consider a differentiable manifold $F$, an open interval $I\subset\mathbb{R}$ and an one-parameter family $\{g_t\}_{t\in I}$ of Riemannian metric on $F$, i.e., a smooth map
$\Lambda: I\times F \rightarrow T_{0,2}(F)$, where $T_{0,2}(F)$ denote the fiber bundle of $2$-covariant tensor
on $F$, such that $\Lambda_{_t}:F \rightarrow T_{0,2}(F)$ is a positive definite metric tensor for all $t\in I$. The product manifold $M=I\times F$ 
can be endowed with the Lorentzian metric given at each point $(t,p)\in I\times F$ by

$$\overline{g}=\beta\pi_{_I}^*(-dt^2)+\pi_{_F}^*(g_t) \quad (\overline{g}=-\beta dt^2+g_t \ \ {\rm in} \ \ {\rm short}),$$

\noindent where $\pi_{_I}$ and $\pi_{_F}$ denote the canonical projection onto $I$ and $F$ respectively, and $\beta \in C^\infty(I \times F)$ is  a positive function. 

We will say that a spacetime is \emph{orthogonal-splitted}, if it is isometric to a spacetime $(M, \overline{g}=-\beta dt^2+g_t )$. Observe that every orthogonal-splitted spacetime is stably causal.

On the other hand, maximal hypersurfaces play a relevant role in General Relativity and Lorentzian Geometry by several reasons. Among them, it should
be emphasized that this class of hypersurfaces has an important
role in the analysis of the Cauchy problem with the purpose of
dealing with simpler constraint equations or solve them,
\cite{Cho}, \cite{Cho2}, \cite[Chap. VI]{Cho3}. Also, in the proof of positivity of gravitational mass, maximal hypersurfaces appear as useful geometric 
objects (see \cite{Sch}).  On the other hand, these hypersurfaces admit another interpretation: they describe the transition between expanding and contracting phases of the universe, in some relevant cases \cite{brill}.
 
A general sumary of several reasons justifying the importance of maximal (and constant mean curvature) hypersurfaces in General Relativity can be found in \cite{M-T}. 

From a mathematical point of view, the maximal hypersurfaces arise naturally  as critical points of the
$n$-dimensional area functional for compactly supported normal
variations. A global and motivating result concerning maximal hypersurfaces is the Calabi-Bernstein theorem.
This theorem asserts that the only complete maximal hypersurfaces in the Lorentz-Minkowski
spacetime, $\mathbb{L}^n$, are the spacelike hyperplanes. This fact was proved by
Calabi \cite{calabi} for $n \leq 4$, and extended to arbitrary dimension by Cheng and Yau \cite{chengyau}. In another spacetimes, the
problem of characterizing maximal hypersurfaces has become an interesting and amusing research topic, and its study can be considered classic.
For instance, in \cite{BF}, Brill and Flaherty replaced Minkowski spacetime with a spatial closed universe, and proved uniqueness at large by assuming ${\mathrm{Ric}}(z,z)>0$ for all timelike vectors $z$ (ubiquitous energy condition). In \cite{M-T}, this energy condition was relaxed by Marsden and Tipler to include, for instance, non-flat vacuum spacetimes. On the other hand, Bartnik, \cite{Bar}, proved very general existence theorems and consequently, he claimed that it would be useful to find new satisfactory uniqueness results. More recently, in \cite{A-R-S1}, Al\'{i}as, Romero and S\'{a}nchez proved new uniqueness results in the class of spatially closed generalized Robertson-Walker spacetimes (which include the spatially closed Robertson-Walker spacetimes), under the Temporal Convergence Condition. Moreover,
in \cite{A-M}, Al\'{\i}as and Montiel improve some of these results and, employing a generalized maximum principle due to Omori \cite{Omori} an Yau \cite{Yau}, give 
a general uniqueness result for the case of complete constant mean curvature spacelike hypersurfaces under certain boundedness of the Ricci curvature. Finally, 
in \cite{RRS}, Romero, Rubio and Salamanca give uniqueness results in the maximal case, for spatially parabolic generalized Robertson-Walker spacetimes, which are open models, whose fiber must have a parabolic universal Riemannian covering (see \cite{RRS} and \cite{RRS2}).

The principal aim of this work is to study  uniqueness of compact maximal hypersurfaces in orthogonal-splitted spacetimes which satisfy certain natural first order conditions, i.e., avoiding curvature assumptions. A characterization of orthogonal-splitted spacetimes which admits a compact spacelike hypersurface 
(i.e., the spacetime is spatially closed) is provided (Prop. \ref{spatiallyclosed}). Then, we will prove that, under certain hypothesis, the compact
maximal hypersurfaces must be level hypersurfaces oat the time function.

In Section 3, we assume certain homogeneity in the expansive or contractive behaviour of the spacetime. More specifically, for each
$v\in T_q F$, $q \in F$, denote by $\tilde{v}$
the lift of $v$ on the integral curve $\alpha_q (s)=(s,q) \in I\times F$. We say that the spacetime $({M},\overline{g})$ is \emph{non-contracting in all directions} 
if $\partial_t \beta \leq 0$ and, $\partial_t \overline{g}(\tilde{v},\tilde{v})\geq 0$, for all $q\in F$ and every $v\in T_qF$. Analogously, we say that the spacetime is
\emph{non-expanding in all directions} if the previous inequalities are reversed. Observe that the spacetime is non-contracting (resp. non-expanding) if and only if 
${\cal L}_{\partial_t}\overline{g}$ semi-definite positive (resp.  semi-definite negative), where ${\cal L}$ denotes the Lie derivative.

If we consider the observer field $U=\frac{1}{\sqrt{\beta}}\partial_t$, the proper time $\tau$ of the observers in $U$ is given by $d\tau=\sqrt{\beta}\,dt$. As a consequence, the assumption $\partial_t\beta\leq 0$ (resp. $\partial_t\beta\geq) 0$  assures that the rate of change (acceleration) of the proper time respect to the time function $t$ is non-increasing (resp. non-decreasing). Observe also that the Lorentzian length $\mid\partial_t\mid=-\sqrt{\beta}$ is non-decreasing (resp. non-increasing) along to the integral curves of $\partial_t$. Observe that all physical interpretations here can be considered locally.

On the other hand, the spacelike assumption $\partial_t \overline{g}(\tilde{v},\tilde{v})\geq 0$ (resp. $\partial_t \overline{g}(\tilde{v},\tilde{v})\leq 0$), guarantees that an observer in $U$ measures non-contraction (resp. non-expansion) in all directions of its physical space.

There are many important examples of orthogonal-splitted spacetimes satisfying theses conditions. For instance, any standard static spacetime 
$M=I\times F$, $\overline{g}=-\beta dt^2+g$, being $\beta:F\longrightarrow (0,\infty)$ a smooth function, obeys clearly our assumptions. 

Several 
generalized Robertson-Walker spacetimes can lie also under this category. Recall that these spacetimes, introduced in \cite{A-R-S1}, are given
by the product manifold of an interval of the real line, $I$ and a Riemannian manifold $(F,g_{_F})$ furnished with Lorentzian  metric 
$\overline{g}=-dt^2+f^2g_{_F}$, where $f: I\longrightarrow (0,\infty)$ is a smooth function. Thus, a GRW spacetime also satisfies the imposed conditions 
whenever the derivative $f'$ is signed. 

As a last example, let us consider the differentiable manifold  $M=I\times F$ endowed with Lorentzian metric
 $\overline{g}=-\beta  (t) dt^2 +f(t) \, g_{_F}$, where $g_{_F}$ is a Riemannian metric on $F$. We can rewrite this metric as $\overline{g}=-d\tau^2 + f(t(\tau)) g_{_F}$. If we assume $f'(t)\geq 0$,  $f''(t)\geq 0$  and $\beta'(t)\leq 0$, the spacetime $M$ favors to model an accelerated expanding universe.
 
\vspace{2mm}

Our first uniqueness result is enunciated as follow (see section 3).

\begin{quote}{\it
Let $(M^n,\overline{g})$, $n\geq 4$n be an orthogonal-splitted  spacetime. Suppose
that $(M,\overline{g})$ is non-contracting (or non-expanding) in all directions. Then each compact maximal hypersurface in $M$ must be a leaf of the orthogonal foliation to vector field $\partial_t$, that is, a level spacelike hypersurface
of the temporal function $t$.}
\end{quote}

Some nice consequences are analized, and several applications to relevant class of spacetimes are provided.

In Section 4, the spacetime has a change in its expansive (or contractive) behaviour. In fact, we will assume that
there exists an instant $t_{_0}$ of the temporal function such that the region of the spacetime $t<t_{_0}$ is non-contracting
in all directions, while the region $t>t_{_0}$ is non-expanding in all directions. Then, the level hypersurface $t=t_{_0}$ 
will be called a \emph{level transition spacelike hypersurface}.
The second main result is presented,

\begin{quote} {\it
Let $(M^n, \overline{g})$, $n\geq 4$, be an orthogonal-splitted  spacetime which has a   level spacelike   transition hypersurface. Then every compact maximal hypersurfaces in $M$
 must be a level spacelike hipersurface.}
\end{quote}

Finally, the last section is devoted to present some applications to Geometric Analysis.

%Our paper is organized as follows. In Section 2 we introduce some definitions and
%the notation to be used for spacelike hypersurfaces in orthogonal-splitted spacetimes. Section 3 is devoted to prove our first principal result and to stablish several applications on some important class of spacetimes.  In section 4 we study the uniqueness problem in the case that the spacetime has a change in its non-contracting or non-expanding behaviour. Finally, in the last section we give some applications to Geometric Analysis.

\vspace{10mm}

\section{Preliminaries} 

Let $(M,\overline{g})$ be an orthogonal-splitted  spacetime.
An immersion of a connected $n$-dimensional manifold $x: S \rightarrow M$ is said to be \emph{spacelike} if the induced metric $g$ on $S$ is
Riemannian. In this case, we will refer to $S$ as a spacelike hypersurface. For each spacelike hypersurface $S$ in $M$ we can
take $N \in \mathfrak{X}^\bot(S)$ as the only globally unitary timelike vector field normal to $S$ in the same time-orientation
of the vector field $-\partial_t$ (i.e., such that $\overline{g}(N,-\partial_t)<0$). From the Lorentzian metric, 
$\overline{g}(N,\partial_t)=-\sqrt{\beta} \, \cosh \theta$, where $\theta$ is the \emph{hyperbolic angle}
between $S$ and $-\partial_t$ at any point of $S$ (see, for instance, \cite[Prop. 5.30]{ON}).
Denoting $\pi_{_I}$ the projection onto the time-coordinate factor, for any spacelike hypersurface $x:S\rightarrow M$, 
we  define the function $\tilde{t}= \pi_{_I}\circ x$ on $S$.

Along this paper, the leaves of the orthogonal foliation to the vector field $\partial_t$ are also called level spacelike hypersurfaces. Note that a spacelike hypersurface $S$ is contained in a level spacelike  hypersurface if and only if the function time $t$ is constant on $S$.

Let $\partial_t^\top:=\partial_t +\overline{g}(N,\partial_t)N$ be the tangential component of $\partial_t$ along $S$. It is easy to see that

\begin{equation} \label{tangentet}
\nabla \tilde{t} =- \frac{1}{\beta} \, \partial_t^\top \, ,
\end{equation} where $\nabla$ denotes here the gradient on $S$.

From the Gauss and Weingarten formulas we have
\begin{equation}
\overline{\nabla}_X Y = \nabla_X Y - g(\mathrm{A}X,Y) N \, ,
\end{equation} for all $X,Y \in \mathfrak{X}(S)$, where $\nabla$ and $\overline{\nabla}$ are the Levi-Civita connections of $S$ and $M$, and $\mathrm{A}$ is the
shape operator associated to $N$,
$$
\mathrm{A}X := -\overline{\nabla}_X N \, .
$$ Recall that the \emph{mean curvature function} relative to $N$ is given by $H:=-(1/n) \mathrm{tr}(\mathrm{A})$. A spacelike hypersurface
is called  \emph{maximal} provided that $H=0$.

On the other hand, a maximal hypersurface is a critical of the  volume variation of  certain family of hypersurfaces (see \cite{Jo}). Indeed, let $x:S \longrightarrow M$ an isometric compact inmersed spacelike hypersurface into the spacetime. We shall speak indistinctly of $S$ and its image if no confusion arise. Let $\omega_{_S}$ be the canonical volume element of $S$ with the Riemannian metric induced from $(M,\overline{g})$. The $n$-volume of $S$ is given by

$$V_S=\int_S\omega_{_S}.$$

\noindent We wish know the variation of $V_S$ and $\omega_{_S}$, when we perform a deformation of the hypersurface $S$.  Let $\xi$ be a differentiable vector field defined  on a neighbourhood of $x(S)$. Then, $\xi$ generates a local transformation one-parameter group $\{\varphi_s\}_{s\in J}$ of local transformations, where $s$ is the canonical parameter and $J\subset\mathbb{R}$ an open interval containing $s=0$. We can define a one-parameter family of isometric inmersions of hypersurfaces $x:S\longrightarrow M$ given by $x_s=\varphi_s\circ x$. Observe that $x_0=x$. The corresponding canonical $n$-volume element for each $x_s$ is noted by $\omega_s$ ($\omega_0\equiv\omega_{_S}$). Then a standart calculation leads to

$$\frac{dV_{x_s}}{ds}(0)=\int_S(\overline{{\rm div}} \xi + \overline{g}(\xi,\vec{H})\, \omega_0.$$

\noindent If $S$ is maximal, the Gauss theorem shows that $\frac{dV_{x_s}}{ds}(0)=0$.

\section{Spacetimes with homogeneity in its expanding-contracting behaviour}

First of all, we study conditions under which an orthogonal splitted spacetime, $(M=I\times F,\overline{g})$ admits
a compact spacelike hypersurface. Observe that if $F$, is compact, then
the spacetime is spatially closed. Moreover, it can be proved the converse, as the following reasoning shows.
Let $S$ be a (connected) compact hypersurface in the spacetime $M$ and $[t_0,t^0]\subset I$ be a closed interval such that $\pi_{_I}(S)=[t_0,t^0]$. Consider the canonical projection projection 
$$\pi_{_F}:S\longrightarrow F\, .$$

Let $TU$ be the fibre boundle of unitary tangent vector on $S$, for all $u\in TU$

$$g(u,u)\leq g_{\pi_{_I}(u)}(d\pi_{_F}(u),d\pi_{_F}(u))=A(t,u)g_{t_0}(d\pi_{_F}(u),d\pi_{_F}(u)).$$

\noindent Thus, the smooth function $A(t, u)$ is defined on the compact $[t_0,t^0]\times TU$. Let $\overline{A}$ be the maximum of $A(t,u)$. so, we have 
$$g(v,v)\leq \overline{A} g_{t_0}(d\pi_{_F}(v),d\pi_{_F}(v)),$$ 

\noindent Hence, up a positive constant factor, it is clear that $d\pi_{_F}$ increases the norms of the tangent vector on $(F,g_{t_0})$. Then, making use of \cite[Lemma 3.3, Chap. 7]{Do}, we obtain that $\pi_{_F}$ is a covering map.

\begin{prop} \label{spatiallyclosed}
i) If an orthogonal-splitted spacetime, $(I\times F, \overline{g})$ admits a compact spacelike hypersurface, then $F$ is compact.

ii) Moreover, if the universal Riemannian covering of $F$ is compact and $S$ is a complete spacelike hypersurface such that,  the function $A(t,u)$, with $t\in \pi_{I}(S)$ and $u\in TU$ is bounded, then $S$ is compact.
\end{prop} Previous result extends Prop. 3.2 in \cite{A-R-S1} properly to a wider family of spacetimes.

\vspace{2mm}

In order to obtain our results, we develop a formula which will be key. This formula is based on the 
computation of the Laplacian of the function $\tilde{t}$ on a maximal hypersurface $S$ inmersed  in a orthogonal-splitted  spacetime $(M, \overline{g})$. 

Let $(S, g)$ be a maximal hypersurface inmersed in the spacetime $M$, where $g$ denotes the induced Riemannian metric on $S$.
Taking 
into account (\ref{tangentet}), a straightforward computation leads to 
\begin{equation}\label{lapla}
\Delta \tilde{t} =  \frac{1}{\beta^2} \, \partial_t^\top (\beta) -\frac{1}{\beta} \, \mathrm{div}(\partial_t^\top )    \, ,
\end{equation} where $\Delta$ and ${\rm div}$ denote the Laplacian and divergence operators  of the inmersion $S$. Now, we can write
$$
\mathrm{div}(\partial_t^\top)= \overline{\mathrm{div}} (\partial_t)+\overline{\mathrm{div}}(\overline{g}(N,\partial_t)N)+ \overline{g}(\overline{\nabla}_N \partial_t^\top, N) \, ,
$$ where $\overline{\mathrm{div}}$  denote the divergence operator of the spacetime. Since the hypersurface
is maximal, we obtain

\begin{equation}\label{div}
\mathrm{div}(\partial_t^\top)= \overline{\mathrm{div}} (\partial_t) + \overline{g}(\overline{\nabla}_N \partial_t, N) \, .
\end{equation}

On the other hand, the term $\overline{\mathrm{div}} (\partial_t)$ has a geometrical interpretation and, in fact, it can be written
in terms of certain volume function. Let $\Omega$ be the canonical Lorentzian volume element (see, for instance, \cite{ON}) of the spacetime. Making use of the Cartan's theorem,
$$
d\circ i_{_{\partial_t}} \Omega = \overline{\mathrm{div}}(\partial_t) \Omega \, .
$$

Consider the level spacelike hypersurface $\{t=t_0\}$, $t_{0}\in I$, and let $\omega_0$ be the canonical Riemannian volume element on $(F, g_{_{t_{_{\small 0}}}})$. Since the spacetime has
an orthogonal splitting, then $\Omega= \sqrt{\beta} dt\wedge \,  \mathrm{vol}_{\mathrm{slice}} \omega_0$, where the function $\mathrm{vol}_{\mathrm{slice}}$ satisfies that $\mathrm{vol}_{\mathrm{slice}}\, \omega_0=\omega_t$, being $\omega_t$ the canonical Riemannian volumen element of $(F,g_{_t})$. Hence, we have 

\begin{equation}\label{div2}
\overline{\mathrm{div}}(\partial_t) = \partial_t \log \mathrm{vol}_{\mathrm{slice}} + \frac{1}{2} \frac{\partial_t \beta}{\beta} \, .
\end{equation}

On the other hand, let us write the normal vector as follows,
$$
N = -\frac{\cosh \theta}{\sqrt{\beta}} \partial_t + N^\perp \, ,
$$ where $\overline{g}(N^\perp,\partial_t)=0$. Hence,
\begin{equation}\label{expresion}
\overline{g}(\overline{\nabla}_N\partial_t,\partial_t)=\frac{\cosh^ 2\theta}{\beta}
\overline{g}(\overline{\nabla}_{\partial_t} \partial_t,\partial_t)+ \frac{\cosh\theta}{\beta}\{\overline{g}(\overline{\nabla}_{\partial_t} \partial_t, N^\perp)+\overline{g}(\overline{\nabla}_{N^\perp}\partial_t,\partial_t)\}+ \overline{g}(\overline{\nabla}_{N^\perp}\partial_t,N^\perp).
\end{equation} 

Let $q\in S$ be such that $N^\perp(p)\not=0$. Take coordinates $({\cal U},(t\equiv x_0, x_1,...,x_n))$ around $p$ in the product manifold $I\times F$, such that $N^\perp=\partial_{x_1}$ on $S\cap {\cal U}$. Therefore, taking into account the equality which determines the Christoffel symbols for a coordinate system, we have

$$\overline{g}(\overline{\nabla}_{\partial_t} \partial_t, N^\perp)=-\frac{1}{2} \,\partial_{x_1} \overline{g}(\partial_t, \partial_t)$$

\noindent and

$$\overline{g}(\overline{\nabla}_{N^\perp}\partial_t,\partial_t)=\frac{1}{2} \, \partial_{x_1} \overline{g}(\partial_t, \partial_t).$$

Furthermore,

$$\overline{g}(\overline{\nabla}_{N^\perp}\partial_t,N^\perp)= \frac{1}{2} \, \partial_t \overline{g}(\partial_{x_1},\partial_{x_1}).$$

Now, using (\ref{div}) and (\ref{div2}) we get

$$
\overline{g}(\overline{\nabla}_N \partial_t, N)= -\frac{\cosh^2 \theta}{2 \beta} \partial_t \beta + \frac{1}{2} \partial_t \overline{g}(\partial_{x_1},\partial_{x_1}) \, .
$$

Finally from (\ref{lapla}) we arrive to
\begin{equation} \label{laplacianotau}
\Delta \tilde{t} = \frac{1}{\beta^2} \partial_t^\top \beta - \frac{1}{\beta} \left[ \partial_t \log \mathrm{vol}_{\mathrm{slice}} -\frac{1}{2} \frac{\sinh^2 \theta}{\beta} \partial_t \beta +       \frac{1}{2} \partial_t \overline{g}(\partial_{x_1},\partial_{x_1}) \right] \, .
\end{equation}

Now, assuming that the spacetime has dimension greater or equal than four, we can consider, on $S$, the pointwise conformal metric to the induced one $\tilde{g}= \beta^{2/(n-2)} g$. In terms of this new metric,
equation (\ref{laplacianotau}) changes to

\begin{equation} \label{laplacianoconformetau}
\tilde{\Delta} \tilde{t} = - \beta^{-n/(n-2)} \left[ \partial_t \log \mathrm{vol}_{\mathrm{slice}} -\frac{1}{2} \frac{\sinh^2 \theta}{\beta} \partial_t \beta +     \frac{1}{2} \partial_t \overline{g}(\partial_{x_1},\partial_{x_1}) \right] \, .
\end{equation}

Observe that the term $\partial_t \log \mathrm{vol}_{\mathrm{slice}}$ can be 
interpreted physically as the increasing or decreasing of the local volume of the physical space for the observer in $U$.

Therefore, if the spacetime is non-contracting or non-expanding in all directions, we have for every point $q\in S$ such that $N^\perp(q)\not=0$,  $\tilde{\Delta} \tilde{t}$ is equally signed.

On the other hand, let $q'\in S$ be with $N^\perp(q')=0$ and assume that there is no an open set ${\cal V}\in S$, with $q'\in {\cal V}$ such that $N^\perp(p)=0$ for all $p\in {\cal V}$. As a consequence, we can take a sequence 
$\{ q_n \}_{n=1}^\infty$ in ${\cal V}$ such that $\lim_{n\rightarrow\infty}q_n=q'$, $N^\perp(q_n)\not=0$ for all $n$ and $\tilde{\Delta} \tilde{t}(q_n)$ has equal sign for all $n$. Thus by continuity $\tilde{\Delta} \tilde{t}(q')$ is equally signed. Finally, if there exists an open set ${\cal V}\subset S$ such that $N^\perp(p)=0$ for all $p\in {\cal V}$, then $\tilde{\Delta} \tilde{t}=0$ on ${\cal V}$.

\vspace{3mm}

We are now in position to state the following result,

\begin{teor} \label{t1}
Let $(M^n,\overline{g})$, $n\geq 4$, be an orthogonal-splitted  spacetime. Assume
that $(M,\overline{g})$ is non-contracting or non-expanding in all directions. Then every compact maximal hypersurfaces in $M$ must be a level spacelike  hypersurface.
\end{teor}

\noindent\emph{Proof.}
From the assumptions, equation (\ref{laplacianoconformetau}) reads as $\tilde{t}$ is a superharmonic or subharmonic function on
a compact Riemannian manifold. Then, $\tilde{t}$ must be constant.
\hfill{$\Box$}

\vspace{2mm}

\begin{rema}\label{vol}
{\rm \textbf{a)} As equation (\ref{laplacianoconformetau}) shows, if such hypersurface $\{t=t_{_0}\}$ exists, then it must be satisfied
$\partial_t \beta (t_{_0})=0$, $\partial_t \log \mathrm{vol}_{\mathrm{slice}} (t_{_0})=0$.
\textbf{b)} A simply calculation shows that
$$\frac{d\omega_s}{ds}(0)=\frac{1}{2} {\rm trace}_{_S}[x^*({\cal L}_\xi\overline{g})] \, \omega_0.$$
Assuming that the level hypersurface $\{t=t_0\}$ is maximal and taking $\xi=\partial_t$, we obtain that
$$\frac{1}{2} {\rm trace}_{_S}[x^*({\cal L}_{\partial_t}\overline{g})]=0.$$ 

\noindent If moreover, $\partial_t$ is a conformal vector field, then $\rho_{\mid_{t=t_0}}=0$, when $\rho$ denotes the conformal  factor.

}
\end{rema}

As the proof of previous result shows, when $\beta \equiv 1$, then the conformal change does not apply and, consequently, we can also
study lower dimensions spacetimes.

\begin{coro}
Let $(I\times F,-dt^2 + g_{_t})$ be an orthogonal-splitted  spacetime, with $\beta\equiv 1$, such that the physical space of the observer in $\partial_t$ non-expands 
(resp. non-contracts) in all directions.
Then the compact maximal hypersurfaces must be a level spacelike hypersurfaces.
\end{coro}

\begin{rema}{\rm Note that when $\beta\equiv 1$ the vector field $\partial_t$ is a observer field, with proper time $t$ and such that the observers in $\partial_t$ are proper time synchronizables. Moreover, each level spacelike hypersuperface $F(t):= \{t\}\times F$ at $t$ constant is a restspace of $\partial_t$, being each event in $F(t)$ simultaneous for all the observer in $\partial_t$.

Let $p,q\in F$ be, the distance between $p_t:=(t,p)$ and $q_t:=(t,q)$ measures by a observer at the instant $t$ is given by $d_t(p_t,q_t)$, where $d_t$ is the riemannian distance in $F(t)$. So, the spacetime is non-contracting (resp. non-expanding) in all directions if and only if $d_t(p_t,q_t)$ is non-decreasing (resp. non-increasing) for all $p,q\in F$.}
\end{rema}

Now, we can provide some consequences of previous results for certain relevant spacetimes. 
The first application that we can consider is relative to multiply warped product spacetimes. Consider
an interval of the real line $I$, $n$ Riemannian manifolds $(F_{_i},g_{_i})$ and $n$ positive funcions $f_{_i}\in C^\infty (I)$
in order to endow the product  $I \times F_{_1} \times \ldots \times F_{_n}$ with the Lorentzian metric
$
-dt^2 + f_{_1}^2\,  g_{_1} + \ldots +f_{_n}^2  \, g_{_n} \, .
$

Multiply warped products are interesting in Cosmology. For instance, it is well-known that the Kasner spacetime can be regarded as a multiply warped product \cite{H}. Also, in the study of $(2 + 1)$ BTZ Black Hole, 
the multiply warped products are considered \cite{Ho}. 

Observe that, if $f_{_i}=f$ for all $i$, then the spacetime is a generalized Robertson-Walker spacetime. 

\begin{coro} \label{staticGRW}
Let $(I\times F_{_1} \times \ldots F_{_m},-dt^2 + \sum_{i=1}^m f_i^2 (t) g_{_{F_i}})$ be a multiply warped product spacetime, such that, for all $i$, the functions $f'_i(t)$ have
a common sign. Then every compact maximal hypersurface in $M$ must be a level spacelike  hypersurface.

\vspace{2mm}

In particular, in every Lorentzian product, $(\mathbb{R}\times F,-dt^2+g_{_F})$ the only compact maximal hypersurfaces are the spacelike level hypersurfaces.
\end{coro}

It can be also treated the case of standard static spacetimes (compare with \cite[Prop. 2]{SeSa}),

\begin{coro} \label{estatico}
In an $n(n\geq 4)$-dimensional standard static spacetime $(I \times F,-h^2 dt^2+g_{_F})$, the only
compact maximal hypersurfaces are the level spacelike hypersurfaces.
\end{coro}

Observe that, in the previous case, every  level spacelike hypersurface $\{t=t_{_0}\}$ is totally geodesic.

\vspace{2mm}

Finally, we can also consider the case of certain Lorentzian twisted products. Recall that the product manifold of an open interval $I$ of the real
line, and a Riemannian manifold $(F,g_{_F})$, endowed with the Lorentzian metric $\overline{g}=-dt^2 + \lambda g_{_F}$, where $\lambda \in C^\infty (I\times F)$ is a  positive function, is a Lorentzian twisted product. In this case, it is well-known that the level hypersurfaces at $t$ constant are totally umbilic (see \cite{PR}).

\begin{coro}
In a Lorentzian twisted product spacetime $(I\times F,-dt^2 + \lambda g_{_F} )$, such that $\partial_t  \lambda \geq 0$ (resp. $\leq 0$)
Then every compact maximal hypersurfaces in $M$ is a level spacelike hypersurfes at $t$ constant.
\end{coro}

\section{Behaviour of the volume function}

In this section, we consider spacetimes whose expanding or contracting behaviour changes its monotonicity. For these purposes,
we introduce the following notion. We say that in an orthogonal-spplitted  spacetime, $(I\times F,-\beta dt^2+g_t)$, a level spacelike hypersurface $t=t_0$ is a  
\emph{level transition hypersurface} if the spacetime is non-contracting in the region $t\leq t_0$ and non-expanding in the region $t\geq t_0$.
Note that if there exist two different spacelike level transition hypersurfaces, then the
portion of the spacetime contained between them has geometrical and physical meaning. In geometrical terms, in this region the
-metric of the- spacetime has no dependance with the $t$-coordinate; in physical terms, the spacetime does not evolve here (it is standard
static in this region). 

This change of behaviour is similar than the spatially closed Friedmann model present (see, for instance, \cite[Ch. 12]{ON}).

It is easy to show that a level spacelike hypersurface has mean curvature $nH=\frac{1}{\sqrt{\beta}} \overline{\mathrm{div}}(\partial_t) -\frac{1}{2\beta^{3/2}}\partial_t \beta$. Then, any level transition hypersurface must be maximal.

Our second main result,

\begin{teor} \label{t2}
Let $(M^n, \overline{g})$, $n\geq 4$, be an orthogonal-splitted  spacetime which has a  level spacelike  transition hypersurface. Then every compact maximal hypersurfaces in $M$
 must be also a  level spacelike hypersurface. 
\end{teor}

\noindent\emph{Proof.}
Let $\{t=t_{_0}\}$ be the level spacelike  transition hypersurface. Consider the function $\tilde{t}-t_{_0}$ in a compact maximal hypersurface. 
From (\ref{laplacianoconformetau}), it is found that
$$
(\tilde{t}-t_{_0}) \tilde{\Delta}(\tilde{t}-t_{_0}) \geq 0 \, .
$$ This fact implies that $\tilde{t}-t_{_0}$ is constant, and so the hypersurface must be a level spacelike hypersurface. To see clear this
last step, consider any function $h \in C^\infty (P)$, where $P$ is a compact Riemannian manifold, such that $h \Delta h \geq 0$. Then, taking into account $\Delta h^2$ and the Divergence Theorem, it is no difficult to see that $h$ must be constant. This applied to the function of this proof, ending it.
\hfill{$\Box$}

\vspace{2mm}

The previous theorem can be applied to a wide family of spacetimes. For instance, let $(F_{_i},g_{_i})$, $1\leq i \leq m$, be a finite
set of compact Riemannian manifolds. Consider the Lorentzian multiply warped product spacetime $(\mathbb{R}\times F_{_1 }\times \ldots \times F_{_m}, 
-dt^2 + \sum_i a_{_i}^2 \exp (-b_{_i}^2 t^2) g_{_i})$, where $a_{_i},b_{_i}$ are not zero constants. As a consequence of the Theorem \ref{t2}, in this spacetime the only compact maximal hypersurface is  $\{t=0\}$. 

\begin{rema} {\rm The concept of transition hypersurface can be given, in a similar way, for a spacelike hypersurface (non-necessarily for a spacelike level hypersurface). Nevertheless a such transition hypersurface may not be maximal, as the following example shows. Consider the spacetime $(I\times \mathbb{S}^1,-dt^2+(-t^2+2t\sin x+3)dx^2)$,
where $I\subset \mathbb{R}$ is the maximal domain of the $t$ coordinate where the spacetime can be defined. Here, the transition 
hypersurface is given by $t=\sin x$.}

\end{rema}

\begin{rema}{\rm If the inequalities are reversed in the definition of transition hypersurface, counterexamples to uniqueness
results can be shown. For example, consider the
 de Sitter spacetime, $(\mathbb{R}\times \mathbb{S}^n,-dt^2+ \cosh^2 (t) g_{_{\mathbb{S}^n}})$. Every rigid motion of the  level spacelike hypersurface $\{t=0\}$ is a maximal hypersurface. Observe that, in this case, $\{t=0\}$ constitutes the  level spacelike transition hypersurface with reversed behaviour.}
\end{rema}

Again, when $\beta\equiv 1$, it can be also treated  the case of dimension $2+1$.

\begin{coro}
Let $(I\times F, -dt^2+g_t)$  be an orthogonal-splitted spacetime which has a level spacelike transition  hypersurface. Then every compact maximal hypersurface in $M$ must be a   level spacelike hipersurface.  
\end{coro}

\section{Applications to Geometric Analysis}

We will begin this section showing that every compact spacelike hypersurface in an
orthogonal-splitted  spacetime is a maximal hypersurface in a suitable
pointwise conformal spacetime. Next, our uniqueness result are applied to study of several Geometric Analysis's problems, which arise naturally from Lorentzian Geometry.

Let $S$ be a spacelike hypersurface isometrically inmersed in an othogonal-spplited spacetime $(M,\overline{g})$. It can be considered the hypersurface $S$ isometrically inmersed in the pointwise conformal spacetime $(M,\tilde{g}=e^{2\alpha}\overline{g})$, where $\alpha \in C^\infty (I\times F)$. The unitary
normal vector field $N$ of $S$ in $(M,\overline{g})$,  
is related with the unitary normal vector field $\tilde{N}$ of $S$ in $(M,\tilde{g}=e^{2\alpha}\overline{g})$  by $\tilde{N}= e^{-\alpha} N$. Observe that if $S$ is spacelike in $(M,\overline{g})$,
then $S$ is also spacelike in $(M,\tilde{g})$. If we take $\left\{E_i \right\}_{i=1}^n$ a local frame field on $S$ in $(M,\overline{g})$, then $\left\{ e^{-\alpha} E_i \right\}_{i=1}^n$, denoted by $\{\tilde{E_i} \}_{i=1}^n$,
is a local frame field on $S$ in $(M,\tilde{g})$ . Denoting by $H$ and $\tilde{H}$ the mean curvature of $S$ in $(M,\overline{g})$
and $(M,\tilde{g})$, respectively, we have
$$
n \tilde{H}=\sum_{i=1}^n\tilde{g}(\tilde{\nabla}_{\tilde{E}_i} \tilde{N}, \tilde{E}_i)= \sum_{i=1}^n g(\tilde{\nabla}_{E_i} e^{-\alpha} N, E_i)=  \sum_{i=1}^ne^{-\alpha} g(\tilde{\nabla}_{E_i} N, E_i ) \, ,
$$ where $\tilde{\nabla}$ is the Levi-Civita connection of $\tilde{g}$. Then
from the previous equation it follows
\begin{equation}\label{curvatura}
e^{\alpha} \tilde{H} = H + \overline{g}(\overline{\nabla} \alpha,N) \, .
\end{equation} On the other hand, we need the following technical result,

\begin{lema}  Let $S$ be a compact spacelike hypersurface inmersed in orthogonal-spplited spacetime $(M,\overline{g})$ and $N$ its unitary normal vector field. Then, for every function $h\in C^\infty(S)$, there exists a function $\overline{\alpha}\in C^{\infty}(M)$ such that 
\begin{equation} \label{NA}
(\nabla \overline{\alpha})|_S^{\bot} = -h N \, .
\end{equation}
\end{lema} 

\textit{Proof.}  The compactness of $S$ allows to take an open interval $J$ such that the geodesics $\gamma_p(s)$ starting at $p\in S$ satisfying $\gamma'_p(s)=N(p)$ are defined for all $p\in S$ and for all $s\in J$. On the other hand, the function $h$ can be extended on a tubular neighbourhood ${\cal U}=\{\gamma_p(s)/ t\in J,\, p\in S\}$  of $S$, being $\overline{h}$ constant along each geodesic $\gamma_p(s)$, with value $h(p)$. Moreover, the vector field  defined in each point of $U$ by $\gamma'_p(s)$, also extends to the normal vector field $N$. Obviously, its flow  on $I\times S$ is given by $\varphi(t,p)=\gamma_p(t)$, being $\varphi$ bijective. Consider the smooth function $\pi_{_J}\circ {\varphi}^{-1}$ Thus, the normal component gradient  $\overline{\nabla}\, \overline{h} (\pi_{_J}\circ {\varphi}^{-1})|_S^\bot = h N$, since $\overline{g}(N,  \overline{\nabla} h ) = 0$.

Now, let $\phi$ a function on $M$ such that $0\leq \phi(p)\leq 1$, for all $p\in M$, satisfying (see Corollary in Section 1.11 of \cite{W}),

\textit{i)} $\phi(p)=1$ if $p\in \{\gamma_t(p)/t\in J',\, p\in S\}$, being $J'\subset J$ an closed interval with $0\in J'$.

\textit{ii)} ${\rm supp}\, \phi\subset U$.

\vspace{2mm}
Therefore, the function $\phi$ allows to extend the function $\pi_{_J}\circ {\varphi}^{-1}$ on all $M$ and to obtain $\overline{\alpha}$.
\hfill{$\Box$}

\begin{lema} In an orthogonal-splitted spacetime every compact spacelike hypersurface is a maximal hypersurface in a suitable pointwise conformal spacetime.
\end{lema}

Let $(M,\overline{g})$ be an orthogonal-splitted spacetime and let $\alpha\in C^ \infty (F)$ a smooth function. Consider the natural extension $\overline{\alpha}=\alpha\circ \pi_{_F}$ on $M$. Let $S$ be a compact maximal hypersurface in $(M,\overline{g})$, then from (\ref{curvatura}) the mean curvature function of $S$ isometrically inmersed in the conformal spacetime $(M,e^{2\overline{\alpha}}\overline{g})$ is given by

$$\tilde{H}=e^ {-\overline{\alpha}}\overline{g}(\overline{\nabla}\overline{\alpha},N) \, .$$

Note that $S$ is the only compact maximal hypersurface in  $(M,\overline{g})$ if and only if $S$ is the only compact spacelike hypersurface in $\overline{\alpha}=\alpha\circ \pi_{_F}$ with mean curvature $\tilde{H}=e^ {-\alpha}\overline{g}(\overline{\nabla}\overline{\alpha},N).$ Making use of this fact, we can obtain some Geometric Analysis's result. 

Let $F$ be a compact Riemannian manifold and let $(I\times F,\overline{g})$ be an $n(\geq 4)$-dimensional orthogonal-splitted spacetime. Let $u:F\longrightarrow I$ a smooth function, such that the entire graph  $\Sigma_u=\{(u(p),p) \, : \, p\in F\}$
is spacelike. Consider $\alpha\in C^ \infty (F)$ and denote us by $N$ the unitary normal vector field on $\Sigma_u$. If we denote $\tilde{H}(u)$ the mean curvature function on the graph $\Sigma_u$ isometrically inmersed in the conformal spacetime $(I\times F,e^ {2\overline{\alpha}}\overline{g})$, with $\overline{\alpha}=\alpha\circ \pi_{_f}$, and we assume that the spacetime $(I\times F,\overline{g})$ is non-contracting or non expanding in all direction, or  it has  a level spacelike transition hypersurface, then every solution to the equation 

$$\tilde{H}(u)=e^ {-\overline{\alpha}}g_{u}(d\pi_{_F}(N), D_u\alpha),$$

\noindent being $D_u$ tthe gradient operator in $(F,g_u)$,  must be constant.

\vspace{2mm}

This last can be specialized to the case in which the ambient Lorentzian manifold is a generalized Robertson-Walker spacetime.

Let $(F,g_{_F})$ be a compact $n$-dimensional Riemannian
manifold and $f : I \longrightarrow \mathbb{R}$ a positive smooth
function. For each $u \in C^{\infty}(F)$ such that $u(F)\subseteq
I$, we can consider its graph $\Sigma_u=\{(u(p),p) \, : \, p\in F\}$
in the Lorentzian warped product $(M=I\times_f
F,\overline{g})$. The graph inherits from $M$ a metric,
represented on $F$ by
\[
g_u=-du^2+f(u)^2g_{_F}.
\]
This metric is Riemannian (i.e. positive definite) if and only if
$u$ satisfies $|Du|<f(u)$ everywhere on $F$, where $Du$ denotes the
gradient of $u$ in $(F,g_{_F})$ and $| Du|^2=g_{_F}(Du,Du)$. Note
that $\tilde{t}(u(p),p)=u(p)$ for any  $p \in F$, and so $\tilde{t}$ and $u$
may be naturally identified on $\Sigma_u$.

When $\Sigma_u$ is spacelike, the unitary normal vector field on
$\Sigma_u$ satisfying $\overline{g}( N,\partial_t)<0$ is
\[
N=\frac{1}{f(u)\sqrt{f(u)^2-\mid D u\mid^2}}\,\left(\,
f(u)^2\partial_t + Du \,\right).
\]

\begin{teor} \label{teonoexistencia} Let $(F,g_{_F})$ be a compact Riemannian manifold, and the smooth functions $\alpha :F\longrightarrow \mathbb{R}$ and 
$f:I\longrightarrow (0,\infty)$, where $I=(a,b)$ is an open interval, satisfying:

i) $f'$ is signed or

ii) $f$ is non-decreasing in the interval $(a,t_0)$ and non-increasing in the interval $(t_0,b).$

\vspace{2mm}
\noindent The only entire  spacelike graphs $\Sigma_u=\{(u(p),p) \, : \, p\in F\}$  in the spacetime $(I\times F, -e^ {2\overline{\alpha}}dt^ 2+e^ {2\overline{\alpha}}g_{_F})$, whose mean curvature $\tilde{H}$ is prescribed by

$$\tilde{H}(u)=\frac{e^{-\overline{\alpha}}}{f(u)\sqrt{f(u)^2-\mid D u\mid^2}}g_{_F}(D\alpha, Du)$$

\noindent must be the constant graph $u=u_0$ such that $f'(u_0)=0$.
\end{teor}

Note that the graph $\Sigma_u$ is spacelike in $(I\times F, -e^ {2\overline{\alpha}}dt^ 2+e^ {2\overline{\alpha}}g_{_F})$ if and only if  it
is so in $(I\times F, -dt^ 2+f(t)^2 \, g_{_F})$. On the other hand, if we conceive  the mean curvature function of a graph $\Sigma_u$ as a known operator $\tilde{H}(u)$ on $C^\infty(F)$, previous theorem can be enunciated in equivalent form,

\begin{teor} Let $(F,g_{_F})$ be a compact Riemannian manifold, $\alpha :M\longrightarrow \mathbb{R}$ a smooth function and $f:I\longrightarrow (0,\infty)$ a smooth function, where $I=(a,b)$ is an open interval,  satisfying:

i) $f'$ is signed or

ii) $f$ is non-decreasing in the interval $(a,t_0)$ and non-increasing in the interval $(t_0,b).$

\vspace{2mm}
\noindent Then, the only solution to the differential equation $$\tilde{H}(u)=\frac{e^ {-\alpha}}{f(u)\sqrt{f(u)^2-\mid D u\mid^2}}g_{_F}(D\alpha, Du),$$ with the restriction $\mid Du\mid<f(u)$ are the constant function $u=u_0$, with $f'(u_0)=0$.

\end{teor}

\end{document}